\documentclass[graybox]{scrartcl}

\usepackage{type1cm}        
\usepackage{makeidx}         
\usepackage{graphicx}        
\usepackage{multicol}        
\usepackage[bottom]{footmisc}
\usepackage{amsthm,amsfonts}

\usepackage{newtxtext}       %
\usepackage[varvw]{newtxmath}       

\makeindex             

\usepackage[utf8]{inputenc}
\newtheorem{theorem}{Theorem}

\begin{document}

\title{Shooting Methods for Fractional Dirichlet-Type Boundary Value Problems 
	of Order $\alpha \in (1,2)$ With Caputo Derivatives}

\author{Kai Diethelm\footnote{\ Faculty of Applied Natural Sciences and Humanities (FANG), 
	Technical University of Applied Sciences Würzburg-Schweinfurt,
	Ignaz-Schön-Str.\ 11,
	97421 Schweinfurt,
	Germany,
	kai.diethelm@thws.de; ORCID 0000-0002-7276-454X
}}

\date{}

\maketitle

\abstract{For the numerical solution of Dirichlet-type boundary value problems associated 
	to nonlinear fractional differential equations
	of order $\alpha \in (1,2)$ that use Caputo derivatives, we suggest to
	employ shooting methods. In particular, we demonstrate that the so-called
	proportional secting technique for selecting the required initial values
	leads to numerical schemes that converge to high accuracy in a very small 
	number of shooting iterations, and we provide
	an explanation of the analytical background for this favourable numerical behaviour.
	}
	


\section{Introduction}

The main goal of this paper is to establish a numerical method for solving boundary
value problems for the nonlinear fractional differential equation
\begin{subequations}
	\label{eq:bvp}
\begin{equation}
	\label{eq:ode}
	D^\alpha y(t) = f(t, y(t))
\end{equation}
of order $\alpha \in (1,2)$ subject to the Dirichlet-type boundary conditions
\begin{equation}
	\label{eq:bc}
	y(0) = b_0, \qquad y(T) = b_1,
\end{equation}
\end{subequations}
on the interval $[0,T]$ with some $T > 0$ and arbitrary given $b_0, b_1 \in \mathbb R$.
In eq.~\eqref{eq:ode}, $D^\alpha$ denotes the Caputo-type differential operator with 
starting point $0$, defined (under our assumption $1 < \alpha < 2$) as
\[
	D^\alpha z(t) = \frac 1 {\Gamma(2 - \alpha)} \frac{\mathrm d^2}{\mathrm d t^2}
					\int_0^t (t - \tau)^{1-\alpha} \left[ 
							z(\tau) - z(0) - \tau z'(0) 
						\right] \, \mathrm d \tau, 
\]
see \cite[Chapter 3]{Diethelm2010a}.

Numerous algorithms for the numerical solution of similar problems have been proposed, 
but they were typically
constructed and analyzed only for the special case that the differential equation
is linear \cite{AlMdallalSyamAnwar2010a,CenHuangXu2018a,CenHuangXuEtAl2020a,%
GraciaStynes2015a,StynesGracia2015a}.
This is a restriction that we do not want to impose. Therefore, our method
will follow a significantly different design principle than the schemes proposed in those
works.

\section{A Shooting Method}

It is well known that fractional differential equations
of order $\alpha \in (1,2)$ that use Caputo derivatives have many properties in common
with classical boundary value problems of order 2.
For this reason, it seems natural to use a numerical approach that has proven to be
efficient and reliable for second-order problems, namely the concept of shooting methods
\cite{Keller2018a}, also for the fractional-order problems under consideration here. 
The principle of a shooting method comprises re-interpreting the given boundary
value problem as an initial value problem whose second initial value $y'(0)$ at the left end point
of the interval of interest needs to be
determined by an iteration process in such a way that the boundary condition at the
right end point of the interval of integration is met. The latter part is achieved
in an approximate sense by a standard numerical solver for fractional initial value problems.
From an abstract point of view, this leads to the following algorithm:
\begin{enumerate}
\item Set $k = 0$ and guess a value $b_{0,k}'$ for $y'(0)$.
\item Numerically compute an approximate solution $y_k$ for the initial value problem 
	comprising the differential equation \eqref{eq:ode} and the initial conditions 
	$y_k(0) = b_0$ and $y_k'(0) = b_{0,k}'$.
\item 
	\begin{enumerate}
	\item If $y_k(T)$ is sufficiently close to the required boundary value $b_1$,
		accept $y_k$ as an approximate solution for the boundary value problem \eqref{eq:bvp}
		and stop.
	\item	Otherwise, determine a new guess $b_{0,k+1}'$ for $y'(0)$, increment $k$ by $1$, 
		and go back to step 2.
	\end{enumerate}
\end{enumerate}
For a complete description and implementation of such an algorithm, one then needs to specifiy
the following components:
\begin{itemize}
\item the selection of the initial guess $b_{0,k}'$ for $y'(0)$,
\item the numerical solver for the initial value problem required in step 2, and
\item the algorithm for computing the next guess $b_{0,k+1}'$ for $y'(0)$ if the approximate solution
	computed in the current iteration was not sufficiently accurate.
\end{itemize}

In particular, we can interpret the last item in the following way: Since all other data is kept fixed,
we can say that the function value $y_k(T)$ of the current approximate solution at the right end point
depends only on the initial guess $b_{0,k}'$, i.e.\ there exists some function $\phi$, say,
such that $y_k(T) = \phi(b_{0,k}')$. To achieve our goal $y_k(T) = b_1$, we thus have to find
a value $b_{0,k}'$ that satisfies the equation 
\begin{equation}
	\label{eq:nonlin}
	\phi(b_{0,k}') = b_1. 
\end{equation}
Therefore, in order to obtain a
fast shooting algorithm, it is essential to solve this (in general nonlinear) equation efficiently.
An evident idea in this context is to use a Newton method; indeed this has been suggested, e.g., 
in \cite{AlMdallalHajji2015a}. However, the computation of the derivatives required for this approach is
not always feasible, and even if the derivatives can be computed, their evaluation may be too costly.
For this reason, we suggest a different approach that is based on the following
result which is a consequence of a recent observation 
about some properties of fractional initial value problems of order $\alpha \in (1,2)$, 
see \cite[Theorem 8]{ChaudharyDiethelmHashemishahrakia}.

\begin{theorem}
	\label{thm:prop}
	Consider the initial value problems
	\begin{equation}
		\label{eq:ivp}
		D^\alpha y_k(t) = f(t, y_k(t)), \qquad
		y_k(0) = y_{0}, \quad
		y_k'(0) = y_{k,1}
		\qquad
		(k = 1, 2)
	\end{equation}
 	for $t \in [0, \tilde T]$
	with $1 < \alpha < 2$ where $f$ is assumed to be continuous and to satisfy a Lipschitz condition in the
	second variable with Lipschitz constant $L$. Assume that $y_{1,0} \le y_{2,0}$ and $y_{1,1} \le y_{2,1}$.
	Then, defining
	\begin{align*}
		\tilde a_*(t) & := \inf_{\tau \in [0,t], y \ne 0} \frac{f(\tau, y + y_1(\tau)) - f(\tau, y_1(\tau))} y
	\intertext{and}
		\tilde a^*(t) & := \sup_{\tau \in [0,t], y \ne 0} \frac{f(\tau, y + y_1(\tau)) - f(\tau, y_1(\tau))} y,
	\end{align*}
	there exists some $T \in (0, \tilde T]$ such that 
	\begin{equation}
		\label{eq:prop}
		(y_{2,1} - y_{1,1}) t E_{\alpha,2}(\tilde a_*(t) t^\alpha)
			\le  y_2(t) - y_1(t) 
			\le  (y_{2,1} - y_{1,1}) t E_{\alpha,2}(\tilde a^*(t) t^\alpha)
	\end{equation}
	for all $t \in [0, T]$.
\end{theorem}

Here, $E_{\alpha, \beta}$ denotes the usual two-parameter Mittag-Leffler function \cite{GorenfloKilbasMainardiEtAl2020a}.
From Theorem \ref{thm:prop}, we can read off two essential messages:
\begin{enumerate}
\item If $T$ is sufficiently small then we can guarantee that the boundary value problem \eqref{eq:bvp} has a unique solution.
\item In view of the special case $t = T$ of eq.~\eqref{eq:prop}, at least in an asymptotic sense, 
	the difference $y_2(T) - y_1(T)$ of the 
	boundary values of two solutions at the right end point of the interval is proportional to the difference 
	$y_{2,1} - y_{1,1} = y_2'(0) - y_1'(0)$
	of the first derivatives of these solutions at the left end point if both solutions themselves coincide 
	at the left end point, 
	i.e.\ if $y_1(0) = y_2(0)$.
\end{enumerate}
The former of these statements lets us conclude that it makes sense to numerically approximate the solution. 
The latter tells us that we can employ the approach introduced in \cite{DiethelmUhlig2023a}
for fractional terminal value problems of order $\alpha \in (0,1)$ and use the so-called 
proportional secting method to determine the next guess for the initial value for $y'(0)$. 
Effectively, this amounts to solving eq.~\eqref{eq:nonlin} with a secant method.

Thus, following the strategy outlined in \cite{DiethelmUhlig2023a},
our algorithm computes the guess $b_{0,k+1}'$ for $y'(0)$ according to the formula
\begin{align}
	\label{eq:nextguess}
	b_{0,k+1}' 
	& = b_{0,k}' + \left( b_1 - y_k(T) \right) \frac{b_{0,k}' - b_{0,k-1}'}{y_k(T) - y_{k-1}(T)} \\
	\nonumber
	& = \lambda_{k+1} b_{0,k}' + ( 1 - \lambda_{k+1} ) b_{0,k-1}'
\end{align}
with
\[
	\lambda_{k+1} = \frac{b_1 - y_{k-1}(T)}{y_k(T) - y_{k-1}(T)}.
\]
Since our approach shown in \eqref{eq:nextguess} 
uses the values $b_{0,k}'$ and $b_{0,k-1}'$, i.e.\ the current approximation for $y'(0)$ and the 
preceding one, to compute the new guess $b_{0,k+1}'$, 
it is necessary to provide two initial guesses to start the algorithm. 
Thus, to complete the description of the shooting method, we need to specify the initial guesses $b_{0,0}'$
and $b_{0,1}'$ for $y'(0)$. 
It has turned out in the similar application discussed in \cite{DiethelmUhlig2023a} that this
approach is very robust with respect to the selection of the starting values and that the number of
iterations that are required until convergence to a given accuracy has been achieved depends
only very weakly on the chosen values of $b_{0,0}'$ and $b_{0,1}'$. Thus we normally use the
particularly simple and computationally cheap selections
\begin{equation}
	\label{eq:initguess1}
	b_{0,0}' = 0
	\quad \mbox{ and } \quad
	b_{0,1}' = \frac 1 T (b_1 - b_0)
\end{equation}
(the latter being the mean value of the first derivative of the exact solution on the given interval).
In the exceptional case when $b_0 = b_1$, this would lead to $b_{0,0}' = b_{0,1}'$ which is not
admissible, and so we then (and only then) replace \eqref{eq:initguess1} by
\begin{equation}
	\label{eq:initguess0}
	b_{0,0}' = 0
	\quad \mbox{ and } \quad
	b_{0,1}' = \begin{cases}
				\phantom{-}1 & \mbox{ if } y_0(T) > b_1, \\
				-1 & \mbox{ if } y_0(T) < b_1.
			\end{cases}
\end{equation}

Finally, in line with the observations of \cite{DiethelmUhlig2023a,FordMorgado2011a}, we suggest to solve
the fractional initial value problems with the selected initial conditions by the
fractional Adams-Bashforth-Moulton method \cite{DiethelmFordFreed2002a,DiethelmFordFreed2004a}
or the second-order fractional backward differentiation formula \cite{Lubich1985a,Lubich1986a}. For the 
sake of simplicity and efficiency, we use a uniform grid for these solvers and the FFT-based implementation
of the chosen scheme \cite{Garrappa2018a,HairerLubichSchlichte1985a} that computes the 
numerical solution on $N$ grid points in only $O(N \log^2 N)$ operations.

\section{Properties of the Method}

The error of the algorithm can be analyzed in a way that is analog to the approach followed
for the terminal value problem solver in \cite{DiethelmUhlig2023a}. 

\begin{theorem}
	\label{thm:err-bound}
	Assume the hypotheses of Theorem \ref{thm:prop}, that $T$ is chosen as indicated
	in that theorem, and
	that the IVP solver approximates the solution of the
	given fractional differential equation \eqref{eq:ode} with 
	an $O(N^{-p})$ convergence order with some constant $p > 0$ 
	for any initial condition. Then 
	\begin{equation}
		\label{eq:err-k-infty}
		\max_{j = 0, 1, 2, \ldots, N} | y(t_j) - y_k(t_j) | = O(|y'(0) - b_{0,k}'|) + O(N^{-p}).
	\end{equation}
\end{theorem}

\begin{proof}
	The proof follows along the lines of \cite[proofs of Theorems 4.1 and 4.2]{DiethelmUhlig2023a}.
\end{proof}

Moreover, as in \cite[Subsection 4.2]{DiethelmUhlig2023a} we can see 
that the algorithm is numerically robust and that it inherits the 
stability properties of the underlying fractional ODE solver. The stability of the
solvers that we have suggested here (Adams-Bashforth-Moulton or fractional BDF)
have been shown to be suitable in \cite{Garrappa2010a,Garrappa2015b,Lubich1985a}.

\section{A Numerical Example}

As a numerical example to demonstrate the performance of our method, we look at the 
differential equation
\begin{subequations}
	\label{eq:ex-bvp}
	\begin{equation}
		\label{eq:fde}
		D^{1.55} y(t)  = \frac 1 {(t+5)^{0.65}}\sin(1.3 \cdot t \cdot y(t)), 
	\end{equation}
a minor variation of a problem discussed in \cite[Example 5.3]{DiethelmUhlig2023a}.
To specify the associated Dirichlet boundary conditions in a suitable way, 
we first note that the uniqueness of the solution of the boundary value problem indeed 
depends on the length of the interval under consideration.
If we choose, for example, the interval $[0, 5.384]$ and compute some solutions
numerically by the BDF2 with a step size $h = 10^{-6}$ then we observe the behaviour shown in 
Figure \ref{fig:ex1}. We have selected two solutions $y_1$ and $y_2$ to the differential equation \eqref{eq:fde}
satisfying the initial conditions $y_1(0) = y_2(0) = 1$, $y_1'(0) = -0.15$ and $y_2'(0) = -0.3$, respectively.
Because the interval is too long to satisfy the conditions required to guarantee the separation of these
solutions, we find that the graphs of $y_1$ and $y_2$ intersect at $t^* \approx 5.384$ with a common
function value $y_1(t^*) = y_2(t^*) = y^* \approx 0.387$. Thus, for such a long interval, the boundary value problem
for the differential equation \ref{eq:fde} with boundary conditions $y(0) = 1$ and $y(t^*) = y^*$ has at 
least two solutions. In such a case, one can usually not predict towards which of the solutions the shooting 
method will converge.

\begin{figure}[htb]
	\centering
	\includegraphics[width=0.8\textwidth]{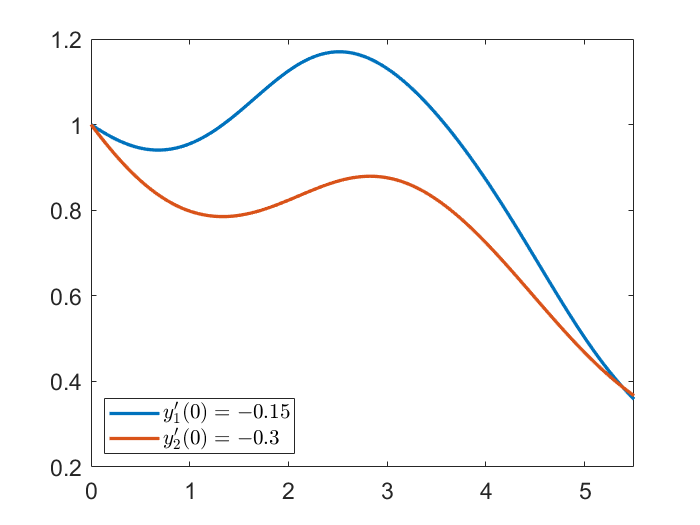}
	\caption{\label{fig:ex1}Graphs of two solutions $y_1$ and $y_2$ to the differential equation \eqref{eq:fde}, both of which satisfy
		$y_k(0) = 1$ but with different values for $y_k'(0)$: The solutions intersect at $t^* \approx 5.384$.}
\end{figure}

Therefore, we need to reduce the length of the interval in order to guarantee the uniqueness of the solution
to the boundary value problem. Following the observations of \cite{ChaudharyDiethelmHashemishahrakia},
the existence of a unique solution can be proved, e.g., if we choose the interval $[0, 2.9]$. Thus, for our numerical
example we combine the differential equation \eqref{eq:fde} with the boundary conditions
	\begin{equation}
		\label{ex:bc1}
		y(0) = 1 
		\quad
		\mbox{ and }
		\quad 
		y(2.9) = 1.145
	\end{equation}
\end{subequations}
that belong to the blue curve in Figure \ref{fig:ex1}.
Using the BDF2 method to solve the initial value problems that arise in this context, 
our algorithm then yields the results shown in Table \ref{tab:ex1}. We can observe that even for a requested accuracy 
of $10^{-10}$, the shooting method converges up to this error in only 7 iterations.

\begin{table}[htb]
	\centering
	\begin{tabular}{c|c|c|c}
		step & requested & no.\ of shooting & actual maximal error \\
		size  & accuracy & iterations & over $[0,2.9]$ \\
		\hline
		$1/100$ & $10^{-4\phantom{0}^{\vphantom{2}}}$ & $5$ & $1.15 \cdot 10^{-5\phantom{0}}$ \\
		$1/100$ & $10^{-8\phantom{0}^{\vphantom{2}}}$ & $6$ & $1.25 \cdot 10^{-7\phantom{0}}$ \\
		$1/200$ & $10^{-8\phantom{0}^{\vphantom{2}}}$ & $6$ & $2.91 \cdot 10^{-8\phantom{0}}$ \\
		$1/400$ & $10^{-8\phantom{0}^{\vphantom{2}}}$ & $6$ & $4.95 \cdot 10^{-9\phantom{0}}$ \\
		$1/400$ & $10^{-10^{\vphantom{2}}}$ & $7$ & $6.76 \cdot 10^{-9\phantom{0}}$ \\
		$1/800$ & $10^{-10^{\vphantom{2}}}$ & $7$ & $7.24 \cdot 10^{-10}$ 
	\end{tabular}
	\caption{\label{tab:ex1}Numerical results for the boundary value problem \eqref{eq:ex-bvp}.}
\end{table}

The computations have been performed with a MATLAB implementation of
our proposed algorithm that can be downloaded from \cite{Diethelm:algorithm}. 

\section{Summary and Conclusion}

In our recent paper \cite{DiethelmUhlig2023a}, we had developed and analyzed the proportional secting scheme for solving 
terminal value problems associated to Caputo-type fractional differential equations of order $\alpha \in (0,1)$.
Here now we change the assumption on the order of the involved differential operators to become $\alpha \in (1,2)$.
Then, the theory of such equations states that a single terminal condition is no longer sufficient to guarantee the
uniqueness of the solution. From our results developed in \cite{ChaudharyDiethelmHashemishahrakia}, we can
see that---under certain conditions on the length of the interval under consideration---the uniqueness can be
recovered if the terminal condition is augmented by an initial condition, thus generating a boundary value problem
with Dirichlet data. It is then natural to ask for an efficient numerical solution method for such equations,
and we have shown here that, by a proper adaptation of the numerical scheme from 
\cite{DiethelmUhlig2023a}, we can create such an efficient algorithm for this higher-order problem with
two supplementary conditions instead of one.

\section*{Software}
The MATLAB source codes of the main algorithm described in this paper
can be freely downloaded from a dedicated repository \cite{Diethelm:algorithm} on the Zenodo 
platform, thus allowing all readers to reproduce the results of our numerical experiments.
The required auxiliary routines are also available from Zenodo; see \cite{DiethelmGarrappaUhlig2023a}.
Our functions were tested in MATLAB R2023a.

\paragraph{Competing Interests} The author has no conflicts of interest to declare that are relevant to the content of this work.


\bibliographystyle{spmpsci}
\bibliography{diethelm-iwota}

\end{document}